\def\strutdepth{\dp\strutbox}
\def \ss{\strut\vadjust{\kern-\strutdepth \sss}}
\def \sss{\vtop to \strutdepth{

\baselineskip\strutdepth\vss\llap{$\diamondsuit\;\;$}\null}}

\newcount\footno
\global\footno=1
\def\comm#1{
\footnote{$^{\the\footno}$}{\sevenrm #1} \advance\footno by 1}

\magnification=\magstep1
\baselineskip16pt
\def\today{\ifcase\month\or
   January\or February\or March\or April\or May\or June\or
    July\or August\or September\or October\or November\or December\fi
\space\number\day, \number\year}
\bigskip
\raggedbottom

\newcount\figno
\global\figno=1
\def\figure#1{\midinsert
\centerline{\epsfbox{#1}}
\centerline{Figure \the\figno}
\global\advance\figno by 1
\endinsert}

\newcount\secno
\newcount\subsecno
\newcount\stno
\global\secno=-1
\outer\def\section#1\par{
     \global\advance\secno by 1
     \global\stno=1
     \global\subsecno=0
     \vskip0pt plus .3\vsize\penalty-250
     \vskip0pt plus-.3\vsize\bigskip\vskip\parskip
     \message{#1}\noindent{\bf{\S\the\secno.\enspace#1}}
     \nobreak\smallskip\noindent}
\outer\def\subsection#1\par{
      \medbreak
      \global\advance\subsecno by 1
      \noindent{\bf#1}
      \nobreak\smallskip}
\outer\def\remark#1\par{
     \medbreak
     \noindent{{\bf Remark   \the\secno.\the\stno.\enspace} {#1}\par
     \global\advance\stno by 1
     \ifdim\lastskip<\medskipamount \removelastskip\penalty55\medskip\fi}}

\outer\def\state#1. #2\par{
     \medbreak
     \noindent{\bf{#1\enspace \the\secno.\the\stno.\enspace}\sl{#2}\par
     \global\advance\stno by 1
     \ifdim\lastskip<\medskipamount \removelastskip\penalty55\medskip\fi}}

\def\tag#1{\edef#1{\the\secno.\the\stno}}
\def\sectag#1{\edef#1{\the\secno.\the\subsecno}}

\newcount\refno
\global\refno=0
\outer\def\ref#1\par{
     \global\advance\refno by 1
     \item{\the\refno.}#1\medskip}

\def\QED{\hfill\rlap{$\sqcup$}$\sqcap$\par\bigskip}
\def\give#1.{\medbreak
             \noindent{\bf#1.}}

\def\g{\gamma} 

\def\more{>}


\def\union{\cup}

\def\{\partial}

\def\int{\text{int}}
\def\invlimit{\smash{\lim\limits_{\raise1pt\hbox{$\longleftarrow$}}}\vphantom{\big(}}
\def\inter{\hskip 1.5pt\raise4pt\hbox{$^\circ$}\kern -1.6ex}
\def\reals{\hbox{I\kern-.13em R}} 
\def\skel(#1,#2){#1^{(#2)}}
\def\hyp {\hbox {\rm {H \kern -2.8ex I}\kern 1.25ex}}
\def\reals {\hbox {\rm {R \kern -2.8ex I}\kern 1.15ex}}
\def\hyp {\hbox {\rm {H \kern -2.8ex I}\kern 1.15ex}}
\def\integers {\hbox {\rm { Z \kern -2.8ex Z}\kern 1.15ex}}
\def\naturals {\hbox {\rm {N \kern -2.8ex I}\kern 1.20ex}}

\def\hyp {\hbox {\rm {H \kern -2.7ex I}\kern 1.25ex}}

\def\HN#1{\parindent0pt \hangafter1 \hangindent35pt \leavevmode%
				 \hbox to 33pt{{\bf #1} \hfil}}
\overfullrule=0pt

\hyphenation{cell-de-com-po-si-tion}
\hyphenation{de-com-po-si-tion}
\hyphenation{di-men-sio-nal}

\font\yoav=cmr8

\input BoxedEPS
\SetTexturesEPSFSpecial
\HideDisplacementBoxes
\centerline {\bf 3-MANIFOLDS WITH IRREDUCIBLE HEEGAARD SPLITTINGS} 
\centerline {\bf OF HIGH GENUS}

\centerline{ Martin Lustig\comm{Supported by Heisenberg-Stipendium from Deutsche 
Forschungsgemeinschaft (DFG).} and Yoav Moriah\comm{Supported by The Fund for Promoting 
Research at the Technion grant 100-053. }}

\bigskip
{\yoav \noindent{\bf Abstract:} Non-isotopic Heegaard splittings of non-minimal genus were 
known previously only for very special 3-manifolds. We show in this paper that they are 
in fact a wide spread phenomenon in 3-manifold theory:  We exhibit a large class of knots
and  manifolds obtained by Dehn surgery on these knots which admit such splittings. Many of the 
manifolds have irreducible Heegaard splittings of arbitrary large genus.  All these splittings are
horizontal and are isotopic, after one stabilization, to a multiple stabilization of certain canonical 
low genus vertical Heegaard splittings.}

\bigskip

\section {Introduction}
\bigskip

Every closed orientable $3$-manifold $M$ has a {\it Heegaard splitting} which is a decomposition  
of $M$ along an orientable surface $\Sigma \subset M$ into two handlebodies
$H_1 , H_2\,$. The genus of this {\it Heegaard surface} $\Sigma$ is called the {\it genus} 
of the splitting. There is a canonical process, called {\it stabilization}, which transforms a
Heegaard splitting of genus $g$ into one of genus $g+1$.  If $M$ is irreducible, then a
Heegaard splitting $M = H_1 \cup_\Sigma H_2$ is {\it irreducible} if it is not obtained
from another splitting of lower genus by stabilization.  A detailed review of these notions
and facts is given below in Section $1$.

The set ${\cal H}(M)$ of all isotopy classes of Heegaard splittings for a given $3$-manifold
$M$ could be determined so far only for a small number of ``simple" manifolds (see
the discussion in Section 1).  Still,   it is known for many manifolds
$M$ that there is more than one isotopy class of {\it  minimal genus}
Heegaard splittings (see [LM1], [LM2]).  For Seifert fibered spaces all irreducible Heegaard
splittings are classified into two types:  They are either {\it vertical} or {\it horizontal} (see
Definition 1.2). There is accumulating evidence that a similar classification might be true
for hyperbolic  manifolds (see [MR], [MS], [CG] and the discussion in Section 1.)

For {\it non-minmal} genus  Heegaard splittings very little is known.  The only manifolds
$M$ for which non-isotopic Heegaard splittings of non-minimal genus have been exhibited
are obtained by surgery on pretzel knots (Casson-Gordon [CG]), or by torus sum of pretzel
link complements  with 2-bridge link complements (Kobayashi [Ko]).  In both cases $M$ is
shown to contain irreducible Heegaard splittings of arbitrary large genus.

In this paper we define {\it vertical} and {\it horizontal} Heegaard splittings in a broad
context, which generalizes the above mentioned earlier notions.  Our results stated
below show that the results of [CG] and [Ko] about high genus irreducible Heegaard
splittings are only the first examples for a phenomenon which is in fact wide spread 
among $3$-manifolds, and which is based on the existence of high genus horizontal Heegaard
surfaces of pairs:

For a general $3$-manifold $M$ and a link $K \subset M$ we introduce the notion of a 
Heegaard splitting of a pair $(M,K)$  which can be  vertical or horizontal (see Section 1.1). 
A  vertical Heegaard splitting of $(M,K)$ will always induce a Heegaard splitting
on all  manifolds obtained by surgery on $K$. However, a horizontal Heegaard splitting 
$\Sigma$ of $(M,K)$,  will induce Heegard splitting only on manifolds obtained by specific
surgeries. Nevertheless we show that horizontal Heegaard splittings are quite common.

Recall that every knot or link  $K \subset S^3$ is isotopic to a $2n$-plat (see Fig. 4.1)  of
length $m$, for some $m,n \in \naturals$, and that every such $2n$-plat can be described by
a family of parameters $a_{i,j} \in \integers$, called {\it twist coefficients} (see Section 4). 
Summing up a well defined subset of these twist coefficients (see Section 4) we  compute a
{\it plat linking number} $a(K) \in \integers$.  To every closed $3$-manifold obtained by 
surgery on a  knot $K$, given as a $2n$-plat, there are two canonically associated 
Heegaard surfaces $\Sigma_{top}$ and $\Sigma_{bot}$ of genus $n$, see Section 5. Let 
$M = K({p \over q})$ denote the closed $3$-manifold obtained from ${p \over q}$-
surgery on $K$.

\tag\Theoremone
\state Theorem.  Let $K$ be knot given as $2n$-plat in  $S^3$, and assume that all twist 
coefficients satisfy $|a_{i,j}| \geq 3$. Then for all  $k \in \integers$, with $|k| \geq 6$,   
 the manifold $K({1 + k \, a(K) \over k})$  has an irreducible Heegaard splitting of genus
$m(n-1)$. Furthermore all these Heegaard splittings are horizontal.

The main tool in this paper is a new combintorial object called {\it trellis},  (see Section 2)
which generalizes the notion of a $2n$-plat and allows us to  present a knot or link by
a family of integer parameters, assembled in a  {\it twist matrix}.  Again, we can compute an
analogous {\it trellis linking number } a(K).   For every knot $K$, carried by a trellis $T$, 
we obtain a {\it trellis Heegaard splitting} of genus $g(T)$ 
for the pair $(S^3, K)$ and for the surgery manifold
$K({1 + k \, a(K) \over k})$. If we consider trellisses with a particular combinatorial feature,
called {\it interior pair of edges}, we can perform flypes at these more  general knots in a
similar way as done by Casson-Gordon in [CG] for  pretzel knots  (see Section
$3$). This allows us to show an analogous result  for a rather large  class of $3$-manifolds:

\tag\Theoremtwo
\state Theorem. Let $T$ be a generalized trellis and let $K = K(A) \subset S^3$ be a  knot 
carried by $T$ with twist matrix $A$. Assume that all coefficients $a_{i,j}$ of $A$ 
satisfy  $|a_{i,j}| \geq 3$ and that there is an interior pair of edges $(e_{i,j}, e_{i,h})$ of 
$T$ with twist coefficients  $|a_{i,j }|, |a_{i,h}| \geq 4$.   Then for all $k, n \in \integers$, 
with   $|k| \geq 6$, the manifolds $K({1 + k \, a(K) \over k})$, have irreducible Heegaard 
splittings $\Sigma(n)$ of arbitrarily large genus $g(T)+2n$, all of which are horizontal.

The above theorems seem, at first sight, to squelch the hope for a natural structure theorem 
concerning the set ${\cal H}(M)$ of {\it all} isotopy classes of Heegaard splittings for
$M$.   However, the following result  perhaps resurrects some of these hopes:

\tag\Theoremthree
\state Theorem. Let $T$ be a generalized trellis and $K \subset S^3$ a knot carried by $T$.
Then for  all $k \in \integers$ the trellis Heegaard splitting of $K({1 + k \, a(K) \over k})$ is 
isotopic, after one stabilization, to a multiple  stabilization of the canonical top Heegaard 
splitting $\Sigma_{top}$ (and also of $\Sigma_{bot}$) defined by $K$. In particular, for
$K$ as in Theorem \Theoremtwo, all of the
splittings $\Sigma(n)$, stabilized once, are stabilizations of a common low genus Heegaard splitting.

\bigskip
Here $\Sigma_{top}$ and $\Sigma_{bot}$ are low genus {\it vertical} Heegaard splittings of
$M = K({1 + k \, a(K) \over k})$ with respect to the core curve $K'$  of the surgery filling torus.  
They generalize the canonical top and bottom splittings for $2n$-plats (see Section 5). 
It has been shown in [LM2] that for sufficiently complicated $2n$-plats, 
$\Sigma_{top}$ and $\Sigma_{bot}$  are typically of minimal genus, and that they are non-isotopic 
in $M$. Examples of arbitrarily high  genus Heegaard splittings which are isotopic after one 
stabilization were found by Sedgwick  (see [Se]). However, it is not known whether the Heegaard 
splittings in those examples are non-isotopic before the stabilization,  nor whether that they are
stabilizations of a common low genus Heegaard splitting.

Haguiwara (see [Ha]) has shown that the canonical top and bottom splittings for $2n$-plats  
become isotopic after at most $2n-1$ stabilizations at each of them.  We give a short proof of his
result in  Section 7, as well as fairly general geometric conditions on the plat which ensure
that less  stabilizations suffice (see Proposition   7.2). 
 
\bigskip
We consider the elements of the set ${\cal H}(M)$ as vertices of a graph in the plane. The
vertices are assembled into horizontal levels according to the genus of the Heegaard splittings.
An edge will connect any two vertices (isotopy classes of Heegaard splittings) if one can be 
obtained from the other by  a single stabilization. The graph  ${\cal H}(M)$ is a  1-ended tree 
(by a well known result of   Reidemeister-Singer), which we call the {\it Heegaard tree} for $M$.  
The results of this paper, as well as all previous results known to us, indicate that ${\cal H}(M)$  
may in general have the following structure:

There is a finite {\it root} part of ${\cal H}(M)$, which contains all irreducible vertical
splitting. Heegaard splittings of the same genus in the root of ${\cal H}(M)$ may 
well need more than one stabilization before they become isotopic, although such phenomenon 
has never been proved so far.  The maximal level of this root part consists  of a single point, 
and from this point there starts an infinite ray moving upward, called the {\it trunk} of  ${\cal H}(M)$. 
At each vertex level of the trunk,  or even of the root part, there may be {\it branches} attached, 
i.e.,  edges which go down into the next  lower level.  Their lower endpoint (a vertex of 
${\cal H}(M)$ not on the trunk)  represents an irreducible horizontal Heegaard splitting. 
In all manifolds known to us these branches all have  length 1.

Since there are only finitely many isotopy classes of Heegaard splittings of the same genus 
(by recent results of Pitts-Rubinstein and Stocking [St]), there are only finitely many such 
branches at each level.  There are examples 
(see [Ko]) where the number of these vertices grows polynomialy, if one moves up the trunk.

\bigskip
\noindent {\bf Aknowledgements:} The authors would like to thank the Technion, the
Frankfurt-Bochum travel grant from the DFG,  and in particular the
Volkswagen-Stiftung's RIP-program at the Mathematisches Forschungsinstitut  Oberwolfach
for their generous help.

\bigskip

\section Heegaard splittings of pairs

\bigskip
In this section we define the basic set up for this paper.  For general definitions and
terminology see  [BZ2], [Ro]  and  [He]. At the end of the section we give a short 
survey about the development of the notions of vertical and horizontal Heegaard splittings.

A compression body  $W$  is a  3-manifold with a prefered boundary component
$\partial_+W$  and is obtained from a collar of  $\partial_+ W$  by attaching 2-handles and
3-handles, so that the connected components of  $\partial_- W$ = $\partial W - \partial_+ W$ are all
distinct from  $S^2$.  The extreme cases, where  $W$  is a handlebody i.e., $\partial_- W = \emptyset$,
or where $W = \partial_+W \times I$, are admitted. Notice that, contrary to the original definition in  [CG], we
require here (as in  [ST]  and  [Sh]) that compression bodies be connected.

A  Heegaard splitting  $(W_1, W_2)$  of a  3-manifold  $M$,
possibly with non-empty boundary, is a decomposition
$M = W_1 \union W_2$, where the  $W_i$  are compression bodies
and  $W_1 \cap W_2 = \partial_+ W_1 = \partial_+ W_2$.  The genus
of the  Heegaard surface  $\Sigma = \partial_+ W_1 = \partial_+ W_2$  is called the
genus of the Heegaard splitting.

A  Heegaard splitting  $(W_1, W_2)$  of a  3-manifold  $M$  is called weakly reducible if there
are disjoint essential disks  $D_1 \subset W_1$  and  $D_2 \subset W_2$. Otherwise it is called
strongly irreducible.  A  Heegaard splitting is called reducible if there are two essential disks
$D_1 \subset W_1$  and  $D_2 \subset W_2$  so that  $\partial D_1 = \partial D_2$; otherwise it will be 
called irreducible.

Given a handlebody $H$, let $D \subset H$ be a collection
of disks so that $H - \inter N(D)$ is a collection of 3-balls. A wave $\omega$ with respect to $D$
is an arc in $\partial H$ so that $\partial \omega$ is contained in a component of
$D$, furthermore $\inter \omega \cap D = \emptyset$, and $\omega$ is not homotopic 
relative $\partial \omega$ into $D$.
\bigskip

\tag\pairdef
\state Definition. \rm Let  $K$  be a knot or a link in a  3-manifold  $M$.  A {\it Heegaard splitting of the
pair   $(M,K)$} is a Heegaard splitting of  $M$  where the Heegaard surface  $\Sigma$ contains  $K$  as
a union of simple closed curves.

\tag\verdef
\state Definition. \rm A  Heegaard splitting for the pair $(M,K)$  is called {\it vertical}  if for each
component of  $K$ there is some properly embedded essential disk in one of the two handlebodies which
is intersected transversally precisely once by $K$ , and it is called {\it horizontal} if $\Sigma - K$ is
incompressible in  $M - K$ (which is the same as saying that it is incompressible in both handlebodies).

Notice that if the genus of a vertical Heegaard surface $\Sigma$ of the pair $(M,K)$ is 
bigger than the number of components of $K$ then $\Sigma - K$ is always compressible. Hence the
vertical and the horizontal case are in this sense opposites extremes of each other.

\bigskip

First examples of a horizontal Heegaard surface $\Sigma$ for a pair $(S^3, K)$ are given by any
incompressible {\it free} Seifert surface $S$ for the links $K \subset S^3$.   Any link in $S^3$ has a free
Seifert surface $S$, i.e., an orientable surface $S \subset S^3$ with $\partial S = K$ such that the
complement of $S$ is a handlebody, and
$\Sigma$ is obtained by simply defining $\Sigma = \partial N(S)$. The Seifert algorithm for obtaining a 
Seifert surface for a knot or link always gives such a free Seifert surface. If  $S$  is incompressible, then
$\Sigma$ will be horizontal. However, in general it is not true that a free Seifert surface will be
incompressible.  In fact, there are knots in  $S^3$  for which any free Seifert surface must be compressible
(see [Ly]).

\bigskip

If $\Sigma$ is a vertical Heegaard surface for the pair  $(M,K)$
then  it gives rise to a  Heegaard splitting for the manifold
$M- \inter N(K)$.  This splitting is obtained by isotopying each component
$K_i$ off $\Sigma$ into the handlebody which contains the disk punctured
once by $K_i$. The handlebodies $H_1$ and $H_2$  are then transformed into
compression bodies $W_1$ and $W_2$  in  $M- \inter N(K)$ which together
determine a Heegaard splitting for  $M- \inter N(K)$.  In particular, since the
components $K_i$ are core cuves of the original handlebodies, this gives a Heegaard
splitting for all closed manifolds obtained by surgery on  $K$, for any surgery value.

If, on the other hand, the Heegaard surface  $\Sigma$ of the pair  $(M,K)$
is not vertical, it will in general not be isotopic in  $M$ to a Heegaard surface
for  $M- \inter N(K)$. The boundary of a neighborhood $N(K_i)$  of each component
$K_i$ of $K$  , $ i = 1, ..., d$, is cut by  $\Sigma$ into two annuli $A^1_i$  and $A^2_i$.
The surface $\Sigma - N(K)$  determines a splitting (not a
Heegaard splitting!) of $M - \inter N(K)$ into two handlebodies $W_1$ and
$W_2$  which are glued along $\Sigma -  \inter N(K) = \partial W_1 - \union A^1_i
= \partial W_2 - \union A^2_i$.  If for each $i$ we glue the two annuli $A^1_i$
and $A^2_i$  together by a multiple Dehn twist along either of them, the
Heegaard surface $\Sigma$ will define a Heegaard splitting of the resulting
manifold.

Let $\beta_i \subset \partial N(K_i)$ be a curve dual to $\partial_i\Sigma = \partial \Sigma
\cap \partial N(K_i)$, i.e., a curve on $\partial N(K_i)$ intersecting $\partial_i \Sigma$ in a
single point. We can choose $\beta_i$ to bound a meridian disk in
$N(K_i)$.   Then, for any integer  $k_i$, glueing the annuli $A^1_i$ and
$A^2_i$ together via a $k_i$-fold Dehn twist is equivalent to performing
$1\over k_i$-surgery on $K_i$ with respect to the basis $(\beta,\partial_i \Sigma)$ 
for $H_1(\partial N(K_i))$, and conversely.  Let us denote by  $\Sigma({1 \over{ \hat k}} )$,
where $ \hat k = (k_1, ..., k_d) $, the manifold obtained by $1\over k_i$-surgery
on each $K_i$.  Hence for  all pairs  $(M, K)$ the surface $\Sigma$ determines a
Heegaard splitting for the pair $(\Sigma({1\over{ \hat k}} ), K)$.

In unpublished work Casson and Gordon proved the following important result
(see [MS] for a proof).  It is formulated here in the terminology introduced above:
\bigskip

\tag\CassonGordon 
\state Theorem. {\rm [Casson-Gordon]}  Let $K \subset S^3$ be a knot and $\Sigma$ a
horizontal Heegaard surface for the pair $(S^3, K)$.  Then for all manifolds
$\Sigma({1\over k})$, with $|k| \geq 6$,  the Heegaard splitting determined by the surface
$\Sigma$ is strongly irreducible.

\bigskip

If $K \subset M$  is a link in a some manifold $M$ and $\Sigma$  is a Heegaard
surface for  $M$, then we say that  $\Sigma$  is {\it vertical (or horizontal) with
respect to} $K$  if  $K$  can be isotoped onto  $\Sigma$  to give a vertical
(or horizontal) Heegaard splitting of the pair  $(M, K)$.  If the reference to
$K$ is non-ambiguous, we sometimes simply say that  $\Sigma$  is vertical
(or horizontal).
\bigskip

We conclude this section by giving some of the history of vertical and
horizontal Heegaard splittings.  None of the following is used in the
sections to come.

\bigskip
{\it Vertical} Heegaard splitting were first defined by Boileau and Otal in
the context of Seifert fibered spaces over $S^2$ with three exceptional fibers.
These are  Heegaard splittings where the handlebodies contain the exceptional 
fibers as cores, i.e., as curves which meet an essential disk in one of the handlebodies
in a single point. It was known already then, by an observation
of Casson and Gordon,  that there were other Heegaard splittings
for these Seifert fibered spaces and that one could isotope an excetional
fiber onto the Heegaard surface also in these cases (see [BO]). The Heegaard
splittings for general Seifert fibered spaces that were described
by Boileau and Zieschang in [BZ1] before the work of Boileau-Otal were by our
definition all vertical with respect to any of the exceptional fibers, 
while the exceptional Heegaard splittings, case (i) of Theorem 1.1 of [BZ1], 
are Heegaard splittings of the pair $(M,f)$ where $f$ is an exceptional fiber.

In unpublished work Casson and Gordon showed that one could find more
examples of horizontal  Heegaard splittings. They showed that some of
the manifolds obtained by surgery on certain generalized pretzel knots
admit irreducible Heegaard splittings, where the core curve of the surgery
torus can be isotoped onto the  Heegaard surface.  The complementary
surface is incompressible to both sides, thus defining a horizontal splitting
of the pair. These knots are all hyperbolic knots, (see [Ka]),  which shows
that horizontal Heegaard splittings are not confined to Seifert fibered spaces.

The view point that these exceptional  Heegaard splittings were infact not
an exotic phenomenon at all was strengthened by a structure theorem for
irreducible Heegaard splittings of negatively curved 3-manifolds,
proved by the second author and Rubinstein (see [MR]). They showed
that given a link in a negatively curved 3-manifold and the full collection
of manifolds obtained by surgery on this link then {\it all}
Heegaard surfaces for "almost all" of these manifolds come from Heegaard
surfaces of the pair  $(M,K)$.

{\it Horizontal} Heegaard splittings were, in the case of general  Seifert fibered
spaces, introduced first by the second author and Jennifer Schultens (see [MS]).
They showed that for orientable Seifert fibered spaces all  Heegaard splittings
are stabilizations of either vertical or  horizontal  Heegaard splittings. Here a
horizontal Heegaard splitting (see Definition 3.1 of [MS]) of a Seifert
fibered space is one which is obtained from a surface fiberation over the
circle of the complement, in the manifold, of a fiber.  Note that in this case
the fiber can be isotoped onto the Heegaard surface and that the Heegaard
surface less a neighborhood of the fiber is incompressible in both handlebodies.
Note also that not all Seifert fibered spaces have such Heegaard splittings.
\bigskip
\bigskip
\eject

\section {\bf  Trellis Heegaard splittings}

\bigskip

Let $T$ be a graph in a vertical plane $P  \subset \reals^3$ which consists
of horizontal and  vertical edges only. Every maximal connected union of
horizontal edges  of  $T$   is called a {\it horizontal line}. The union of
two adjacent horizontal lines and all  vertical edges spanned between them is
called a {\it horizontal layer}. If $T$  has  $m$ horizontal  layers and contains  $n$ vertical edges 
in each of them,  arranged in \lq\lq brick like fashion" as in Fig. 2.1, it is called a {\it standard trellis} 
of size $(m, n)$.   Its regular  neighborhood in $\reals^3$  is a handlebody $H_1 = N(T)$ of 
genus $g(T) = m(n-1)$, embedded in the standard way in $ S^3$, which we identify with the one-point
compactification of  $\reals^3$.

\vskip25pt
\hskip80pt\BoxedEPSF{Fig.2.1.eps scaled 700}
\bigskip
\centerline {Fig. 2.1}
\bigskip

For any integer $(m \times n)$-matrix  $A = (a_{i.j})$ we define a knot or link
$K = K(A)$ contained in the boundary of the handlebody $H_1$  and winding
around the trellis  $T$ as in  Fig. 2.2. There each configuration as in Fig. 2.3, occuring at the  
$j$-th twist box, counted from  the left, of the $i$-th layer, counted from the top, indicates
$a_{i,j}$  half twists. We call  $a_{i,j}$ the {\it twist coefficients} and $A$ the {\it twist matrix}.  
We always use $P$ as projection plane for $K = K(A)$.

\eject

Note that the long horizontal strings are on the back of the trellis.

\vskip12pt

\hskip60pt\BoxedEPSF{Fig.2.2.eps scaled 720}

\bigskip

\centerline{Fig. 2.2}

\vskip20pt

\hskip130pt\BoxedEPSF{Fig.2.3.eps scaled 680}
\medskip
\centerline{Fig. $2.3$}

\bigskip

Whenever a  trellis  $T \subset S^3$  and a knot or link  $K \subset S^3$ are
given, and  $K$  is isotopic to  $K(A)$  for some twist matrix   $A$  as
above, then we say that  $K$  is {\it carried by}  $T$ {\it with twist matrix} $A$.

As the complement  $H_2 = S^3 - \inter N(T)$  of  $H_1$ is also a handlebody,
the pair $(H_1,H_2)$  defines a  Heegaard splitting of the pair  $(S^3, K)$,
which we  call the {\it trellis Heegaard splitting}. We refer to $H_1$ as the
{\it inner} handlebody  and to $H_2$ as the {\it outer} handlebody of this
splitting. As in the last section we denote the surface which is their common boundary  
by $\Sigma$.  The plane $P$ cuts $\Sigma$ into two connected components which we refer
to as the {\it front} and the {\it back}.

Notice that  $K$  bounds a possibly non-orientable surface  $S$  in  $H_1$,
defined by replacing every vertex of  $T$  by a small disc in  $P$  and every
edge of  $T$  by a twisted band attached  to those disks.

\bigskip

\tag\innerhandlebody
\state Lemma.  The surface $\Sigma - K$  is incompressible in the inner handlebody $H_1$.

\give Proof. The handlebody  $H_1$  admits a structure of an orientable
$I$-bundle over the  surface  $S$. Hence  $\Sigma - K$  is isotopic in $\partial H_1$ to the 
induced orientable $\partial I$-bundle over  $S$, which is the orientable double cover of $S$  in
case  $S$  is  unorientable, or, if  $S$   is orientable, it is the disjoint union of two
copies of $S$.  In particular the fundamental group of (a component of) $\Sigma - K$ is mapped
injectively to $\pi_1(S) = \pi_1(H_1)$.

\QED

\tag\outerhandlebodythree
\state Lemma.  If $n \geq 3$  and if all twist coefficients satisfy   $|a_{i,j}| \geq 3$   
then $\Sigma - K$  is incompressible in the outer handlebody  $H_2$.

\give Proof. Notice that the projection plane  $P$ cuts the handlebody $H_1$
through the middle. Let ${\cal D} = \{ D_{i,j}\}$, for $1\leq i \leq m$ and  $1\leq j \leq n-1$, be 
the set of disks given by those connected components of   $P \cap H_2$  which are compact.
The complement of ${\cal D}$  in $H_2$  is a 3-ball.

We now want to remove all inessential intersections of $K$ with $D$ by an isotopy
of $K$ on $\Sigma$ (\lq\lq tightening $K$ with respect to $D$"). Such inessential
intersections occur only at the top or the bottom horizontal line of $T$. It occurs 
exactly if some $a_{1,j}$  is negative and $a_{1,j+1}$ is positive, or if $a_{m,j}$ 
is positive and $a_{m,j+1}$ is negative. Hence by isotopying some of the top  and 
some of the bottom arcs of $K - P$ from the front to the back of $\Sigma$ we eliminate
all of the inessential intersections. Notice that our assumption  $|a_{i,j} | \geq 3$  implies 
that after these isotopies each {\it vertical column} of $H_1$, i.e., the neighborhood of a
vertical edge of the underlying  trellis, has at least one small horizontal arc of $K$ on
the front of  $\Sigma$ which connects the  two adjacent disks  $D_{i,j}$ and
$D_{i, j+1}$, and another  such arc on the back.
  
Let $\gamma$  be a loop in  $\Sigma - K$ which is contractible in  $H_2$ and
transverse to ${\cal D}$.  Hence, after tightening $\gamma$ with respect to ${\cal D}$,  the 
loop $\gamma$  either misses  ${\cal D}$  or else it  contains a  wave with respect to $D$
(see Section 1). 

It follows from our assumption $| a_{i,j} | \geq 3$  that the connected components $\Sigma_i$ 
of $\Sigma - \inter N(K - \cal {D})$  do not contain essential loops.  Therefore the curve  
$\gamma$ must meet ${\cal D}$ and hence it must contain a wave $\gamma_0$.

This wave $\gamma_0$ is a properly embedded arc in some connected component
$\Sigma_i$  of  $\Sigma - \inter N(K - \cal {D})$. Its endpoints are on some $\partial D_{i,j}$,
(more precisely, on the parallel copies of $\partial D_{i,j}$ on $\partial N(D_{i,j})$),
and $\gamma_0$ must join two distinct connected  components of the intersection $D_{i,j} \cap
\Sigma_i$  from the same side of $D_{i,j}$. Thus we have reduced our goal, showing that there is 
no  essential loop in $\Sigma - K$ which bounds a disc in $H_2$, to showing  the following:
\smallskip
\noindent {\it Claim:}  For each connected component  $\Sigma_i$  of
$\Sigma - \inter N(K - \cal {D})$ the intersection of $\Sigma_i$ with any  $D_{i,j}$ is
either empty, or consists of precisely one arc, or consists of precisely  two
arcs along which  $\Sigma_i$ meets  $D_{i,j}$ from opposite sides of $P$.

To prove this claim we devide the complementary components of  $K \cup \cal
{D}$  in $\Sigma$  into finitely many classes, according to their position on
$\Sigma$ as pictured  in Fig. 2.4.

\vskip15pt
\hskip10pt\BoxedEPSF{Fig.2.4.eps scaled 660}
\medskip
\centerline {Fig. 2.4}
\bigskip

Those complementary components which are just small horizontal strips on the
front or the back   of a vertical column of $H_1$ satisfy the claim, as they
meet $\cal D$  in precisely two  arcs,   which belong to distinct  $D_{i,j}\,$ unless
the column is an outermost one. For the outermost columns  the horizontal strips
run around from the front to the back and hit the same disk $D_{i,j}$ twice,  but
from opposite sides of $P$, so they also satisfy the claim.

Next we consider the class of complementary regions  $\Sigma_i$ which are located on the front of 
$\Sigma$, and are in one to one correspondence with the valence-3-vertices of  $T$  other than
those on  the top or on the bottom horizontal line. Each such \lq\lq triangular shaped" $\Sigma_i$ 
can meet at most  four  disks $D_{i,j}$, and these are  all distinct, unless $\Sigma_i$  is
outermost on its horizontal layer. In the latter case we notice that the
assumption  $n \geq 3$  implies that the triangular region can not be
outermost simultaneously to the right and to the left.  Hence at most two of
the four intersection arcs may belong to  the same disk $D_{i,j}$, but then
$\Sigma_i$   meets that disk $D_{i,j}$ from opposite sides of $P$, which proves
the claim for this class of regions as well.

A third class of complementary components $\Sigma_i$ arises on the back of $\Sigma$. Each 
such $\Sigma_i$ contains in its boundary one of the subarcs of $K$ which have been isotoped
from the front to the back in the tightening process of $K$, and $\Sigma_i$ meets at most three  distinct 
$D_{i,j}$. If $\Sigma_i$ is outermost on its  horizontal layer two of the $D_{i,j}$ will agree, but 
then they are met  by $\Sigma_i$ from opposite sides of $P$ .

It remains to check the last class, consisting of \lq\lq long horizontal" complementary
regions, one on the top front of the first layer, one on the bottom front of the last
layer, and two regions on the back of each layer .  However, it is easy to check that 
each of those regions meets any non-outermost  (in its horizontal level) disc $D_{i,j}$ 
in at most  one arc, while the outermost discs $D_{i,j}$  could be met by some  regions possibly 
twice, but if that happens then they are met  from opposite sides of $P$.
This proves the Claim and hence the lemma.

\QED

In  what follows we will admit more general knots  $K = K(A)$ than the ones considered so
far:  We start with a standard trellis $T$ of size $m \times n$ and remove any number of
vertical edges or horizontal edges, with the following restrictions: There are at least
three vertical edges in each layer, there is  only one horizontal line in each horizontal level,
there are no edges of valence  one, and the trellis is connected.  The resulting graph
$T'$ will be called  a {\it generalized trellis}.  As before, a knot or link carried by
the generalized trellis defines a twist matrix $A$, which is an $(m \times n)$-matrix
$A$ with integer coefficients, except that we use the convention that
we set $a_{i,j} = \infty$ for those entries of $A$ which correspond to the vertical edge
$e_{i,j}$ of $T$ that were deleted when passing over to $T'$.

Conversely, given such a matrix $A$, the knot or link $K = K(A)$ is built in the
neighborhood of the deleted edges on the local model  used for standard  trellisses at the top
and at the bottom horizontal lines, so that all the terminology   and all the  basic facts for
standard trellisses extend to the case of a generalized trellis $T$ as well. We define the {\it
genus} $g(T)$ of $T$ to be the genus of the  handlebody $H_1 = N(T)$.  

The proof for the incompressibility of  $\Sigma - K$ in the
inner handlebody $H_1$ (as in Lemma \innerhandlebody) carries over word by word to
generalized trellises. In order to prove the incompressiblity of $\Sigma - K$ in the outer handlebody
$H_2$ we need to make the  following adjustments in the proof of  Lemma \outerhandlebodythree:

\item {(a)}  In the tightening process of $K$ with respect to $D$, we may
need to isotope additional arcs from the front to the back. These additional arcs will  
occur at the top or the bottom of the deleted vertical edges.

\item {(b)} We will now have to consider regions which replace the
triangular shaped components of $\Sigma - \inter N(K - D)$  on the front of $\Sigma$
for a standard trellis, but are more complicated. These regions may now have more \lq\lq sides": 
They do not necessarily  correspond any more  to single vertices on interior horizontal
lines of $T$, but rather to segments on such lines. These segments contain only vertices 
which bound vertical edges from above or only vertical edges from below, and are maximal 
with respect to this property. It is easy to see that these new regions still satisfy the {\it Claim} in the 
proof of Lemma \outerhandlebodythree , so that  the proof of this lemma carries over directly
to generalized  trellises $T$.

\smallskip
\noindent We summarize the results of this section with the following:

\medskip
\tag\horizontalthree
\state Proposition.  Let  $K = K(A)$ be a knot or link  carried by a
generalized trellis  $T$  with twist matrix $A$ that has coefficients $a_{i,j} \in
\integers \union \{ \infty \}$.  If all twist coefficients satisfy $|a_{i,j} | \geq 3$  
then the trellis Heegaard splitting of the pair $(S^3, K)$   associated to $T$ is horizontal.

\QED

The proofs of Lemma 2.2 and Proposition 2.3  show that the condition 
$|a_{i,j} | \geq 3$ is by no means a necessary condition for both statements. For example, a local
necessary condition is that not both of  $a_{i,j}$ and $a_{i,k}$ be $0$
for $j\ne k$. However, it seems difficult, at this stage, to give precise necessary and
sufficient conditions.

\bigskip
\tag\Deflink
\state Definition. \rm For any knot $K$ carried by a generalized trellis $T$ we define the 
{\it trellis linking number} $a(K)$ as the element of $H_1(S^3 - N(K))$ determined by a
component of $\partial( \Sigma -\inter N(K))$, where $\Sigma = \partial N(T)$. 

The trellis linking number $a(K)$ can be computed as follows: Choose an orientation for $K$. 
Let $A$ is the twist matrix of $K$.  Define $A_0$ to be the set of all twist coefficients  
$a_{i,j} \in A$
with the property that the two oriented strings of the knot $K$ cross through the corresponding 
twist box of the trellis projection in the same vertical direction. Notice that  in this case the local linking 
of the knot with a parallel curve on the surface $\Sigma$ is twice the twist coefficient $a_{i,j}\,$.
If the orientations of the strings are opposite then the linking number is $0$. Notice also that the strings 
of $K$ on the back of the trellis do not contribute to the local linking. Hence 
$a(K) = 2\Sigma\{a_{i,j} |  a_{i,j} \in A_0\}$. In particular the boundary slope on $\partial N(K)$ 
determined by $\partial( \Sigma -\inter N(K))$, expressed in the usual meridian/longitude coordinates
of $H_1(\partial N(K))$, is $a(K) \over 1$.  
It follows that $\Sigma({1 \over k}) =  K({{1 + k\, a(K)} \over k})$. 
\eject

\noindent{\section {Flypes}

\bigskip

Let $T$ be a generalized trellis. We say that two adjacent vertical edges $e_{i,j}, e_{i,k}$ 
in the $i$-th interior horizontal layer, with $i \ne 1, m$, is an  {\it interior pair of edges}  if  
$e_{i,j}$ and $e_{i,k}$ are not outermost, and if the segments of the two horizontal lines 
bounded by the vertices of $e_{i,j}$ and $e_{i,k}$  satisfy the following properties:
\noindent
 \item{(a)} There are no vertical edges in the $(i-1)$-th or in the $(i+1)$-th
horizontal  layer which have endpoints on either of the two segments.

 \item{(b)} There are  two vertical edges in the  $(i-1)$-th and two in the
$(i+1)$-th  horizontal layer which have endpoints separating the two segments
from the endpoints  of all other vertical edges in the $i$-th layer. This is illustrated  in Fig. 3.1.

\vskip15pt
\hskip100pt\BoxedEPSF{Fig.3.1.eps scaled 640}
\medskip
\centerline{Fig. 3.1}
\medskip

A{ \it flype} at the interior pair of edges $(e_{i,j},  e_{i,k})$   is an ambient
isotopy of  $K$ which is obtained as follows: Consider a box in $S^3$ which intersects 
$K$ in exactly the two subarcs of  $K$ winding around the edge  $e_{i,j}\,$, the two subarcs
winding around $e_{i,k}\,$, and in the two horizontal subarcs on the front of $\Sigma$
connecting the top of $e_{i,j}$ to the top of $e_{i,k}\,$, and the bottom of $e_{i,j}$ to the bottom
of $e_{i,j}$ respectively.  (see Fig. 3.2.)

\vskip15pt
\hskip80pt\BoxedEPSF{Fig.3.2.eps scaled 550}
\medskip
\centerline{ Fig. 3.2}
\bigskip   

A flype will flip the box  by 180 degrees about a horizontal axis leaving all parts of the  knot outside  the
box fixed.  This operation changes the projection of $K$ in $P$ by adding a crossing on the left and 
a crossing on the right side of the box. These crossings have opposite signs.

The projection of $K$ obtained after a flype is carried by a new trellis. It differs
from $T$   in that there is a new vertical edge on the left side of $e_{i,j}$ and
another new one on the right side of $e_{i,k}$, one with twist coefficient $1$
and the other one with $-1$.  The  flype will  be called  {\it positive} if the
coefficient of the right \lq\lq new" edge is positive. A positive/negative flype
iterated $\pm \, n $ times  will be called an $\pm \, n${\it -flype}, (see Fig. 3.3).
When the interior pair of edges in which the $n$-flype is performed is specified before we will 
denote the image of $K$ after the  $n$-flype by $K(n)$ and the new trellis with the new 
$2n$ vertical edges by $T(n)$. Similarly we will denote $N(T(n))$ by $H_1(n)$ and
$\partial H_1(n)$ by $\Sigma(n)$.

As before, the inner handlebody $H_1(n)$  is cut by the projection plane
$P$ through  the middle, and the compact components of the intersection of
$H_2(n) =  S^3 -H_1(n)$  with $P$  give a collection $D(n)$ of disks, which cut
$H_2(n)$ open to give a $3$-ball.  The disk collection $D(n)$ consists precisely of 
the disks $D_{i,j}$ defined as in the last section for $T$, and, for each flype, an additional 
two disks, one on the left of $e_{i,j}$, and one on the right of $e_{i,k}$.

Our next goal is to show that the surface $\Sigma(n) - K(n)$ is incompressible in $H_2(n)$.
As in the last section, we first have to tighten $K(n)$ with respect to the disk
system $D(n)$. This is  done by moving some arcs from the front to the back
part of $\Sigma(n)$,  as explained in the last  section for generalized trellises.
In this tightening procedure we first isotope those arcs from the  front of $\Sigma(n)$ 
to the back which already had to be moved in order to tighten $K$  with respect to 
$D$.  The only place where $K(n)$ may not be tight, after these \lq\lq old" tightening
isotopies, are horizontal arcs with one endpoint on the top or on the bottom of
the vertical column corresponding to the edges $e_{i,j}$ or $e_{i,k}$. This is
because all new left vertical edges have the same sign for their twist coefficient, and  
similarly for all new right edges (with opposite sign).  There are various cases according
to the sign of the  twist coefficients $a_{i,j}$ and $a_{i,k}$, and the sign of
the flype, and they will be discussed in the proof of Lemma 3.2 below.

If we try to  show the incompressibility of $\Sigma(n) - K(n)$ in the outer
handlebody as before we will quickly run into a problem, as it will turn out that often the 
disk system $D(n)$ decomposes $\Sigma(n) - K(n)$ into subsurfaces and some of
them  {\it do indeed} contain a wave.  Thus we first need to generalize our method:

\medskip For any knot or link  $K \subset \Sigma$ and a disk system  $D$ which cuts
$H_2$ into  one (or more) $3$-balls let us consider, as before, a decomposition  of
$\Sigma$ into  subsurfaces $\Sigma_i$ which are simply connected and which have 
boundary on $K \union \partial D$.  We require as before that  $\Sigma_i$ meets $K$  only
in proper  subarcs of $\partial \Sigma_i$, but contrary to the above we allow the possibility
that $\Sigma_i$ contains some properly embedded arcs from $\partial D$. In other words,
the decomposition considered here arises from the connected components  of $\Sigma - (K
\cup D)$ by gluing together some of these components along segments of 
$\partial D$.

Let $\gamma$ be a  path which is properly embedded in $\Sigma_i$ and transverse to $D$,
with boundary  points on two distinct components of
$\partial \Sigma_i \cap D$.  Notice that, as $\Sigma_i$  is simply connected, up to a
homotopy of $(\gamma, \partial \gamma)$ in $(\Sigma, \partial \Sigma -K)$ there are  only
finitely many such paths. We read off the word corresponding to the intersections of 
$\inter\gamma$ with the disks from $D$, and freely reduce it to get the {\it interior word} 
$w (\inter \gamma)$. Let $w(\gamma)$ be the analogously defined word, but with the two
extra intersections of $\gamma$ with $D$ at the two boundary points of $\gamma$.  These
words are   invarint modulo free reduction, with respect to relative homotopy of $\gamma$.
As
$\Sigma_i$ is simply connected these words only depend on the two components of 
$\partial \Sigma_i \cap D$ which contain the endpoints of $\gamma$.

\bigskip

\tag\Testlemma
\state Lemma.  If for each such $\gamma$ the freely reduced words $\omega(\gamma)$ and
$\omega(\inter\gamma)$ satisfy the equation
$$\hbox{length}(w(\gamma)) = \hbox{length}(w(\inter\gamma) + 2 ,$$  then $\Sigma - K$
is incompressible in the outer handlebody $H_2$.

\bigskip
\give Proof.  Every loop $\rho$ in $\Sigma - K$, after being made transverse and tight with
respect to $D$, decomposes into arcs $\gamma_i$ as above, which are concatenated  along
their boundary points: $\rho = \gamma_1 \gamma_2 \ldots \gamma_q$.  By assumption 
any letter which correspond to one of these boundary points, say $\g_i \cap \g_{i+1}$ (with
$i$ understood mod $q$), does not  cancel with either of the adjacent reduced interior words  
$w(\inter \gamma_i)$ or $w(\inter\gamma_{i+1})$  (or against the first letter coming from
the next  arc $\gamma_i$, in case the interior word is empty).

 Hence the whole loop reads off a reduced word which is non-trivial, and thus it can not be
contractible in $H_2$.

\QED
\bigskip Notice that if no $\Sigma_i$ contains any properly embedded arc from $\partial D$
then Lemma \Testlemma\/ coincides with the old criterion that no $\Sigma_i$ may contain a
wave.
\bigskip 

This lemma will be applied below in a particular situation, which we want to
spell out explicitly.  It is easy to see that in this situation the hypotheses of Lemma
\Testlemma \/ are satisfied for the regions $\bar \Sigma_i$ defined below.
\bigskip

\tag\Almostnowaves
\state  Lemma.  Assume that one has a decomposition of  $\Sigma$  along
$K \cup D$ as before. Assume also that for each wave $\gamma_0$ in any of 
the components $\Sigma_i$  the two adjacent regions $\Sigma_j , \Sigma_k$ of
$\Sigma_i $ which contain the endpoints of  $\gamma_0$  satisfy the following conditions:
\hskip 100pt
\item{(a)}  $\Sigma_j$ and $\Sigma_k$ do not contain waves,
\item{(b)}  $\Sigma_j$ and $\Sigma_k$  do not meet any of the curves  $\partial D_{i,j} 
\in \partial D$ from the same side, and
\item{(c)}  the union $\bar \Sigma_i$ of $\Sigma_i$ with all adjacent regions which
contain an endpoint of a wave in 
\hskip 20pt 
\indent $\Sigma_i$ is simply connected. \hfil\break
\noindent Then  $K$ is incompressible in the outer handlebody  $H_2$.

\QED

\medskip

We are now ready to prove:
\bigskip

\tag\lemmatwo
\state Lemma.  Let $K = K(A)$ be a knot or link carried by a generalized trellis
$T$. Assume that all coefficients $a_{i,j}$ of $A$  which are different from
$\infty$  satisfy  $|a_{i,j}| \geq 3$,  and that there is an interior pair of edges
$(e_{i,j}, e_{i,k})$ of  $T$  with $|a_{i,j }|, |a_{i,k}| \geq 4$. Then for 
any $n$-flype at $(e_{i,j}, e_{i,k} )$ the surface $\Sigma(n) - K(n)$ is  
incompressible in $H_2(n)$ .

\give Proof. As the flype involves only a local part of the trellis and the knot or link
carried by it,  we can use the fact shown in the proof of Lemma \outerhandlebodythree
\hskip 3pt  that those components $\Sigma_i(n)$ of $\Sigma(n) - (K(n) \union
\partial D(n))$ which  have not been changed by the flype do not contain a
wave. Thus it suffices  if we investigate those \lq\lq new" components $\Sigma_i(n)$ which
intersect the flype  box defined above. We will have to distinguish various cases
according to the sign of the  flype number $n$ and of the twist coefficients
$a_{i,j}$ and $a_{i,k}$.  In each of  these cases there will be the following types
of \lq\lq new" complementary components of $K(n) \union D(n)$ in $\Sigma(n)$:

\noindent (a)  Small \lq\lq horizontal" strips on the front or on the back of  the
vertical columns corresponding to $e_{i,j}$ and $e_{i,k}$.

\noindent (b) Two long horizontal regions on the front, which bound all of the new
disks: One of the long regions bounds from above and the other long region from below.

\noindent (c) Two similarly long horizontal regions on the back.

\noindent (d) Regions on the back which bound  one of the arcs of $K(n)$ moved
to the  back by our tightening isotopies above, and which bound precisely three
disjoint disks from $D(n)$.

\vskip15pt
\hskip5pt\BoxedEPSF{Fig.3.3.eps scaled 700}
\smallskip
\centerline{Fig. 3.3}
\bigskip 
We now distinguish the following  three  cases, pictured in Fig. 3.3.
All other possibilities  can be treated similarly to one of them, by the two
mirror-symmetries of the given set up.

\noindent I. \hskip 15pt  $n  \leq 1$ and $a_{i,j} \geq 4$ and $a_{i,k} \leq -4$

\noindent II. \hskip 11pt $n  \leq 1$ and $a_{i,j} \leq -4$ and $a_{i,k} \geq 4$

\noindent III. \hskip 8pt $n  \leq 1$ and $a_{i,j} \geq 4$ and $a_{i,k} \geq 4$

In cases  I  and  II  there is precisely one region of type (d), and in case (3)
there is none.  I any case, such regions never contain a wave.
Clearly the regions of type (a) or (d) never contain a wave.   In
case I  we check from Fig. 3.3  that none of the four long horizontal regions of type
(b)  or  (c) contains a wave. In case  II  there are two such regions with precisely
one wave each, namely the bottom region of type (b), and the top region of type
(c).  The other two long horizontal regions do not contain waves.  Similarly,
in case  III  there are two  long horizontal regions which contain a wave:  The
bottom region of type (b), and the top region of type (c).

Observe that in each case the two adacent  regions
which contain the endpoints of the wave $\gamma$ are always of type (a), and the two 
never belong to the same vertical column. It is easy to check that the conditions of 
Lemma \Almostnowaves \ are satisfied, which proves that the surface $\Sigma(n) - K(n)$
is  incompressible in $H_2(n)$.

\QED

\tag\Secondthm
\state Theorem. Let $T$ be a generalized trellis and let $K = K(A) \subset S^3$ be a 
knot or link carried  by $T$ with twist matrix $A$. Assume that all coefficients $a_{i,j}$ of $A$  which
are different from $\infty$ satisfy  $|a_{i,j}| \geq 3$ and that there is an interior pair of edges $(e_{i,j},
e_{i,h})$ of  $T$ with twist coefficients $|a_{i,j }|, |a_{i,h }| \geq 4$.   Then for all $n \in
\integers$ the trellis $T(n)$ obtained from an $n$-flype at this edge pair defines  a trellis Heegaard 
splitting for the pair $(S^3, K)$ which is horizontal and of genus $g(T(n)) = g(T) + 2n$. \hfil\break 
\indent In particular, if $K$ is a knot, then for all the manifolds $K({1 + k \, a(K) \over k})$  
with  $|k| \geq 6$ this induces a strongly irreducible 
Heegaard splitting of genus $g(T) + 2n$, for all $n \in \integers$.
\bigskip

\give Proof.  For all $n \in \integers$ the surface
$\Sigma(n) - K(n)$ is incompressible in the inner handlebody $H_1$  by Lemma
\innerhandlebody, and it is incompressible in the outer handlebody $H_2$ by
Lemma \lemmatwo.  Hence the trellis splitting defined by  $T(n)$ is horizontal,
and as a consequence of  Casson-Gordon's result, stated in Theorem \CassonGordon \/, this 
gives a strongly irreducible  Heegaard splitting of genus $g(T) + 2n$ for the surgery manifolds 
$K({1 + k \, a(K) \over k})$, with $|k| \geq 6$.  

\QED

\noindent{\bf Proof of Theorem 0.2.}  The theorem follows directly from Theorem \Secondthm.

\QED

\bigskip

\bigskip

\noindent{\section {Horizontal Heegard splittings for  knots in plat projections}

\bigskip

\bigskip

In this section we apply the tools developed in the previous two sections
for knots or links carried by a trellis, to knots or links given as plats  (see [BZ] and Fig. 4.1).

\vskip15pt
\hskip60pt\BoxedEPSF{Fig.4.1.eps scaled 650}
\bigskip
\centerline {Fig. 4.1}
\bigskip
A $2n$-plat projection as above, determines  a  $(m \times n)$-{\it prematrix} $ \hat A$ with 
integer twist coefficients $a_{i,j}$. A  $(m \times n)$-{\it prematrix} is a $(m \times n)$-\lq\lq matrix" 
where the odd numbered rows have only $n - 1$ entries instead of $n$. Precisely,  for $i$ odd we have  
$1 \leq j \leq n - 1$, while for $i$ even we have   $1 \leq j \leq n $.

A prematrix $ \hat A$  determines, in a canonical way, a matrix  $A$ by
defining $a_{i,n} = 0$ for all odd indices $i$. We will say that $A$ is obtained from
$\hat A$ by {\it 0-filling}. A first observation is the following:
\bigskip

\tag\plat
\state Lemma. Every knot or link $K$ given as  $2n$-plat with twist prematrix $\hat A$ 
is isotopic to the knot or link $K(A)$ carried by a standard trellis of size $(m,n)$,
with twist matrix $A$ obtained by $0$-filling from $\hat A$.
\bigskip

\give Proof.
 For every odd layer of the plat projection one takes the left-most vertical subarc $k$ of  
$K$  and moves it by an ambient isotopy along the back of the plat projection until it
is in a position right of the former right-most vertical subarc arc in this layer.  This isotopy
creates two long horizontal subarcs on the back, connecting the top end point in the
old position of the arc $k$ to the top of the new position, and similarly at the bottom of
$k$. We now interpret the two right-most parallel vertical strings of this new projection 
of $K$ as the $n$-th twist  box of this layer (with twist coefficient equal to 0), and observe 
that this gives a knot or link $K(A)$  precisely as claimed (see Fig. 4.2 ).

\QED
 
\vskip10pt
\hskip60pt\BoxedEPSF{Fig.4.2.eps scaled 800}}
\smallskip
\centerline{Fig. 4.2}
\bigskip
\bigskip

\tag\plattrellis
\state Proposition. Let $K$ be a knot or link in a $2n$-plat projection, let $\hat A$ be
the associated twist prematrix, and assume that all twist coefficients satisfy $|a_{i,j}| \geq 3$. 
Then the trellis Heegaard splitting of the pair $(S^3, K(A))$ is horizontal, where the 
twist matrix $A$  is obtained by $0$-filling from $\hat A$.

\give Proof. Let $T$ be the standard $(m, n)$-trellis which carries $K(A)$, and
let $\Sigma$ be the associated trellis Heegaard surface. By Lemma 2.1 the subsurface
$\Sigma - K$ is incompressible in $H_1 = N(T)$. Thus it remains to show that
$\Sigma - K$  is incompressible in $H_2 = S^3 - \inter H_1$.  The proof uses
the same technique as that of Lemma 2.2 .

From the assumption that all twist coefficients of the prematrix $\hat A$ associated to the
$2n$-plat $K$ satisfy $|a_{i,j}| \geq 3$ it follows that the twist matrix $A$ for the trellis  $T$
satisfies the same condition, except that $a_{i,n} = 0$ whenever $i$ is odd.  As in the
proof of Lemma 2.2  we consider the decomposition of $\Sigma$ into subsurfaces $\Sigma_i$ 
by cutting along $K \union D$, where $D$ is the disk system in $H_2$ considered there.  It is
shown there that,  if {\it all} twist coefficients of $A$ satisfy $|a_{i,j}| \geq 3$, then none of the
subsurfaces $\Sigma_i$ contains a wave.  Hence it suffices to consider only those subsurfaces
which meet the right-most vertical column of an odd horizontal layer.

It is easy to see that there are exactly three such complementary regions, and that the two
of them which intersect this vertical column only on the front do not contain a wave. 
However, the third one does contain waves on the back of $\Sigma$.  For each disk
$D_{i,j} \in D$ in this layer, except for the right-most, there is a wave. It starts at the top of
$D_{i,j}$, runs horizontally to the right, then down over the right-most vertical
column, and then horizontally back to the bottom of $D_{i,j}$.  Its two  endpoints are in  
different connected components of $\partial D_{i,j} - K$.  A picture is given in Fig. 4.3.

\vskip25pt

\BoxedEPSF{Fig.4.3.eps scaled 760}
\bigskip
\centerline{ Fig. 4.3}
\bigskip

Notice that the waves pointed out  above are the only waves in this region.  Hence we can easily 
verify that the hypotheses of Lemma \Almostnowaves\ are satisfied which implies the incompressibility of 
 $\Sigma - K$ in $H_2\,$.

\QED

Given a knot $K$ is in a $2n$-plat projection we can compute, as in Definition \Deflink\/, 
$a(K)$ with respect to the standard trellis $T$ given by  Lemma $4.1$ \/. It is the linking number of 
a boundary  component of the corresponding 
surface $\Sigma - \inter N(K)$ with $K$. In this case we will call $a(K)$
the {\it plat linking number}.

\bigskip

\noindent{\bf Proof of Theorem 0.1.}  The surface $\Sigma - N(K)$ is incompressible in the inner
handlebody $H_1$  by Lemma \innerhandlebody, and it is incompressible in the outer handlebody 
$H_2$, by Propostion 4.2.  Hence the trellis splitting defined by  $T$ is horizontal, and as a consequence 
of  Theorem \CassonGordon \ this gives a strongly irreducible Heegaard splitting of genus $m(n - 1)$ 
for all surgery manifolds $K({1 + k \, a(K) \over k})$, with $|k| \geq 6$.

\QED
\eject
\noindent \section {Canonical Heegaard splittings}

\bigskip

The goal of this section is to extend the notion of canonical top and bottom Heegaard splittings
which are defined for plat projection of knots or links $K$, reviewed below, to arbitrary
projections of $K$.

\bigskip

\tag\candef
\state Definition. \rm Let $K$ be a knot or link given as $2n$-plat in $\reals^3 \subset
S^3 = \reals^3 \union \{\infty\}$.  Let $\tau_1, ..., \tau_{n - 1}$ be a system of pairwise 
disjoint horizontal arcs which connect adjacent local maxima (the {\it top bridges}) of $K$ 
(see Fig. 5.1). One defines the compression body  $W_1$  to be the union of a collar of  
$\partial N(K)$ in $S^3 - \inter N(K)$,  and of a regular neighborhood of the  
$\tau_1, ..., \tau_{n - 1}\,$.  The handlebody  $H_2 = W_2$  is defined as complement  
$(S^3 - \inter{N(K)}) - W_1\,$, and together they define the {\it canonical top Heegaard splitting} 
of $S^3 - \inter{N(K)}$ associated to the given $2n$-plat projection of $K$.
Analogously, if we replace the arcs $\tau_i$ by similar arcs  $\rho_1, ..., \rho_{n - 1}$
connecting adjacent local minima of $K$ we obtain the  {\it canonical bottom Heegaard splitting}
of $S^3 - \inter{N(K)}$.

\vskip10pt

\hskip60pt\BoxedEPSF{Fig.5.1.eps scaled640}
\bigskip
\centerline{ Fig. 5.1}
\bigskip

\bigskip
\tag\canrem
\state Remark. \rm  As both canonical  Heegaard splittings are obtained by adding
tunnels to a regular neighborhood of $K$ both splittings are (after pushing $K$ on $\partial
N(K)$) vertical splittings for the pair $(S^3, K)$.

\bigskip

Let  $K \subset \reals^3 \subset S^3 = \reals^3 \union \{\infty\}$ be a knot or a link, were we  think 
of $K$ as a specific embedding, rather than its isotopy class, and assume that with respect to the 
standard hight function in $\reals^3$ there are finitely many local maxima of $K$ occuring on
small subarcs $\mu_1, \ldots, \mu_r$ of $K$.  We allow the degenerate case that such a $\mu_i$ 
is a horizontal arc, and we assume that the arcs $\mu_i$ are labeled so that $i \more j $ implies that 
the hight of $\mu_i$  is bigger or equal to the hight of $\mu_j$. We consider a (large) horizontal 
disk $\Delta$ above $K$ and connect every $\mu_i$ by a monotonically ascending arc $\nu_i$ 
to $\Delta$ (see Fig. 5.2 below).  We require that all $\nu_i$ are  pairwise disjoint and do not 
meet $K$ other than at their lowest point (the {\it initial} point). 

We first want to show that the complement $H_2$ of the handlebody $H_1 = N(K
\cup \{\nu_1, \ldots, \nu_r\} \cup \Delta) \subset S^3$ is also a handlebody:  This can be seen
by contracting each $\nu_i$ while moving $\mu_i$ monotonically upward, until it hits
$\Delta$. The result is the disk  $\Delta$  with  $2r$ strands attached on its bottom side which 
descend monotonically until they reach a local minimum of $K$.  These strands are braided, 
but there is an ambient isotopy of $\Delta$ which moves their endpoints around so that the 
braid becomes trivial. This moves $H_1$ into a standard position in $S^3$, and hence the 
complement $H_2$ is also a handlebody.

Next we want to show that the isotopy class of  $\Sigma = \partial H_1=  \partial
H_2$ in $\reals^3 - \inter N(K)$ does not depend on the particular choice of the 
arcs $\nu_i\,$:  For the top arc  $\nu_1$ this is clear, as there is only one isotopy class 
of monotonically ascending arcs.  For the second top most arc $\nu_2$  there is
more than one possible isotopy class, but it is easy to see that the various choices can be 
obtained from each other by sliding the terminal point of $\nu_2$ over $\Delta \cup \nu_1$.  
Similarily we slide  $\nu_3$ over $\Delta \cup \nu_1 \cup \nu_2$  to get all possible isotopy 
classes for $\nu_3$, and so on.  The isotopy class of  $\Sigma$ in $\reals - \inter N(K)$ is 
not changed during these moves, which proves our claim.  This justifies the following:
\bigskip
\tag\canonicalsplit
\state Definition. \rm The above Heegaard splitting  $(H_1 -  \inter N(K), H_2)$  of 
$S^3 - \inter N(K)$  is called the {\it canonical top Heegaard splitting} of the given 
knot or link $K$ and is denoted by $\Sigma_{top}$.  Similarly, if we invert the hight function
i.e., replacing maxima by minima and making the other corresponding changes, we define the {\it
canonical  bottom Heegaard splitting} of the given knot or link $K$,  denoted by $\Sigma_{bot}$.
Notice that these Heegaard splittings depend on the actual embedding of the curve $K$ in
$S^3 = \reals \cup \{\infty\}$ and not just on its ambient isotopy class.

We now want to change the viewpoint slightly:  Suppose $\{\mu_s,\ldots, \mu_t\}$ is a subset
of $\{\mu_1,\ldots, \mu_r\}$ contained in the same horizontal plane $Q$. Consider any monotonically
ascending subarc $k$ of $K$ which crosses $Q$ transversely and connects it to some $\mu_d$ on a
strictly higher level. We isotope all of the arcs $\nu_s, \ldots, \nu_t$ by sliding their terminal point down
along $\nu_d$ and then along $k$ (keeping them throughout pairwise disjoint and their interiors
disjoint from $K$) until they become horizontal arcs $\nu'_s, \ldots, \nu'_t$ contained in the plane 
$Q$.  Furthermore we allow iterative slides of any of the $\nu'_i$, within $Q -K$, over any other $\nu'_j$. 
Again, these slides do not change the isotopy class of the resulting Heegaard splitting.  

In this way we obtain an alternative description of the top canonical Heegaard splitting, defining 
$H_1$  as neighborhood of $K$ and of a system of horizontal and vertical arcs. In particular this 
shows the following:

\bigskip
\tag\samename
\state Remark.  For the special case of a $2n$-plat $K$ the above defined canonical top 
Heegaard splitting $\Sigma_{top}$ coincides, up to an isotopy in $S^3 - \inter N(K)$, with the 
canonical top Heegaard splitting  associated with a $2n$-plat.

We consider now the case of a knot or link $K$ carried by a generalized (!) trellis
$T$ and compare its canonical top Heegaard splitting $\Sigma_{top}$ to the top Heegaard 
splitting $\Sigma_{top}(n)$ of the knot projection $K(n)$ obtained from $K$ by an $n$-flype as defined
in Section 3. 

Consider the $2n$ local maxima arcs $\mu_i$ of $K(n)$ on the same hight level which
are generated by the $n$-flype. They are contained in some horizontal plane $Q$ and are
connected by vertical arcs $\nu_{i_{k}}\,$, $k = 1, \ldots , 2n\,$, to the disk $\Delta$. 
Let  $\mu_{i_{0}}$ be the horizontal local maximum arc between the two vertical strands 
on the interior pair $(e_{i,j}, e_{i,k})$ at which the flype is performed, and let $\nu_{i_{0}}$
be the corresponding vertical arc. Now slide these arcs $\nu_{i_{k}}$, $k = 0, \ldots , 2n\,$, 
as described above so that they become pairwise disjoint horizontal arcs on $Q$ 
(as indicated in  Fig. 5.2 below).

Now undo the flype, and obtain a system of arcs $\nu_{i_{k}}^*$ with endpoints on $K = K(0)$ 
which looks as follows:  There is a horizontal arc $\nu_{i_{0}}^*$ (corresponding to $\nu_{i_{0}}$), 
together with  $n$ vertical arcs on the left and  $n$ vertical  arcs on the right of the interior  pair 
$(e_{i,j}, e_{i,k})$.  The tunnel $\nu_{i_{0}}^*$ either connects the two vertical strands winding 
around $e_{i,j}$, or else those winding around $e_{i,k}$. In the first case (the second is similar) we 
can slide one  endpoint of each of the vertical arcs $\nu_{i_k}^*$ on the left of $(e_{i,j}, e_{i,k})$
over a subarc of $K$  winding around  the  edge $e_{i,j}$ and over the arc $\nu_{i_0}^*$  to
transform it into a trivial tunnel  (see Fig 5.3).

\vskip15pt

\hskip30pt\BoxedEPSF{Fig.5.2.eps scaled640}
\medskip
\centerline {Fig. 5.2}

\vskip20pt

\hskip70pt\BoxedEPSF{Fig.5.3.eps scaled670}
\medskip
\centerline {Fig. 5.3}
\bigskip

Then we can slide the left endpoint of the arc
$\nu_{i_0}^*$ up along $K$ and some $\mu_j$, then through the disk $\Delta$ and over some 
of the $\nu_l$ of higher index, and finally back down on some other vertical arc $\mu_h$ and 
a subarc of $K$ so that it reaches a position where it is a small horizontal arc which connects the 
two vertical strands of $K$ which wind around $e_{i,k}\,$. We then do the same arcs slides on 
the right of $(e_{i,j}, e_{i,k})$ as we did before on the left and hence also transform the other 
$n$ arcs $\nu_{i_k}^*$ into trivial tunnels, thus proving the claim.

This transforms all the $2n$ vertical arcs $\nu_i$ into trivial tunnels: Each $\nu_i$ is a small arc with
boundary on $K$ which bounds together  with a small subarc of $K$ a small disk in $H_2$ and hence
meets a cocore disk in $H_1$  transverse to $\nu_i$ precisely in one point. We obtain:

\bigskip
\tag\Flypestab
\state Proposition.  If the knot or link $K$ is carried by a generalized trellis $T$, and if $K(n)$ is
obtained from $K$ by an n-fold flype at some interior pair, then the canonical
top Heegaard splitting  of  $K(n)$  arises from that of $K$ by $2n$ stabilizations.

\QED 

We finish this section by considering the change of the canonical top Heegaard
splitting induced by adding vertical tunnels:

\bigskip

\tag\horizontalarcs
\state Lemma. Let $K \subset \reals^3$  be a knot or link, and let  $\{\nu_i, \nu'_j\}$ be
a set of horizontal or vertical arcs which determine the
canonical top Heegaard splitting.  Let $\{\alpha_j\}$ be a set of horizontal arcs
with endpoints on $K$  which are pairwise disjoint and meet  $K\cup \{\nu_i, \nu'_j \}$
only in their endpoints.  Then the resulting surface $\partial N(K \cup \{\nu_i\}
\cup \{\alpha_j\})$ is a Heegaard surface of $S^3 - \inter N(K)$,  and it arises from
multiple stabilization of the canonical top Heegaard surface $\Sigma_{top}$ for $K$.

 \give Proof.  We first bring the horizontal arcs  $\nu'_j$ into a
monotonically ascending position $\nu_j$, by succesively sliding one of their endpoints over some of the
other $\nu'_k\,$, some of the $\mu_i\,$, and over one of the $\nu_k$ until it reaches $\Delta$.  This can be
done without changing the position of the $\alpha_j$. Next we contract the arcs $\nu_i$ by sliding the
$\mu_i$ up until they hit $\Delta$  (as described in the begining of the section). We then  move the
$\alpha_i$ iteratively up, starting always with the top-most one, until they  reach $\Delta$. There they
form a collection of  trivial arcs with endpoints on $\Delta$, which proves the claim. 

\QED

\bigskip

\noindent{\section {Stabilizing horizontal Heegard splittings}

\bigskip  Given a  Heegaard splitting of a $3$-manifold $M$, we can obtain a new
Heegaard  splitting by adding a pair of 1-handles, one to each handlebody, so that
their cocore  disks intersect in a single point.  We will call this operation a {\it stabilization}.  
It is well known that any two Heegaard surfaces $\Sigma$  and $\Sigma'$  of $M$ 
become isotopic after a sufficiently large number of stabilizations on both Heegaard surfaces.   
If $q \geq 0$  or less such stabilizations on either surface suffice for such an isotopy, we will 
say  that $\Sigma$   is  $q$-isotopic to $\Sigma'$.  In general it is difficult to determine the
minimal  possible such $q$; an upper bound depending linearily on the genus of
the two surfaces has been given recently in  [RS].
\smallskip

\noindent {\bf Note:} Throughout this section we will always assume that $K$ is a knot.
\smallskip

Part (a) of the following statement seems to be known; for completeness we include a proof.

\bigskip

\tag\Propositionstab
\state Lemma.
\item {(a)}  Let $K \subset M$ be a knot. Every Heegaard surface $\Sigma$  of a pair
$(M,K)$  is  1-isotopic in $M$ to a vertical  Heegaard surface of $(M, K)$ and in 
particular to a Heegaard surface $\Sigma^*$ of $M - \inter N(K)$.
\item{(b)}  Let $K \in S^3$ be a knot and $\Sigma$ a Heegaard surface of the pair
$(S^3, K)$. Let  $K'$ the core of the surgery filling in $\Sigma({1 \over k})$,  and
$\Sigma_k$ the Heegaard surface of the pair $(\Sigma({1 \over k}), K')$ defined
by $\Sigma$ as in Section 1.  Then $\Sigma^*$ and $\Sigma_k^*$, defined as in part (a) 
for  $M = S^3$ and $M = \Sigma({1 \over k})$ respectively, are isotopic to each other in 
$S^3 - \inter N(K) = M - \inter N(K')$, for all $k \in \integers$.
\bigskip

\bigskip

\give Proof.  (a)  Let $H_1$ and $H_2$ be the two handlebodies of the Heegaard
splitting of $M$ given by $\Sigma$. We choose a regular neighborhood $N(K)\subset M$ 
and a meridional disk $D$ for $N(K)$. Consider the arc $\tau = \partial D \cap H_1$ .
It is properly embedded in $H_1$  and has endpoints on different sides of $K$  in  
$\Sigma \cap N(K)$, see Fig. 6.1 below.  We drill out a small neighborhood $N(\tau)$ 
of $\tau$ from $H_1$ and add it to $H_2$, to obtain two new handlebodies $H_1^*, H_2^*$ with
common boundary   $\partial (H_1 - \inter N(\tau)) = \partial (H_2 \union N(\tau))$. This defines a new
Heegaard splitting of the pair $(M, K)$  which is of genus one higher than the original one.  It is a
stabilization of $\Sigma$, as a cocore disk $D'' \subset N(\tau) \subset H_2^*$ meets the disk $D' = D \cap
H_1^*$ exactly once. The disk $D'$ meets  $K$  precisely in one point.  Hence the  new Heegaard
surface is a vertical Heegaard surface with respect to $K$.  If we isotope the new surface
slightly off $K$ into $H_2^*$ we obtain the desired Heegaard surface $\Sigma^*$ of
$M - \inter N(K)$.

\noindent(b) We can canonically identify $S^3 - \inter N(K)$ with $\Sigma({1\over k}) - \inter
N(K')$.   After we push the surface $\Sigma$ off $K$  into $H_2 \subset S^3 - \inter N(K)$ it is
canonically identified with the surface $\Sigma_k$ pushed off $N(K')$. This from from the fact
that the closed manifold $\Sigma({1 \over k})$ can be obtained by gluing $H_1$ to $H_2$ along 
their boundaries, but with the gluing  map  from $S^3$ modified through $k$-fold Dehn twist at $K$,
see Section 1. 

As explained above in (a), the surface $\Sigma^*$ is obtained from $\Sigma$  using an arc
$\tau \subset S^3 - \inter N(K) = \Sigma({1\over k}) - \inter N(K')$ with endpoints on $\Sigma$, 
and  $\Sigma^*_k$ is obtained similarly  from $\Sigma_k$  by an analogous arc 
$\tau_k \subset \Sigma({1\over k}) - \inter N(K') = S^3 - \inter N(K)$.  These two arcs differ 
essentially in that  $\tau_k$ runs once around $\partial N(K)$, as does $\tau$, but in addition 
$\tau_k$ runs  $k$ times parallel to $K$. However, we can define an isotopy between $\Sigma^*$ 
and  $\Sigma^*_k$ by sliding one  \lq\lq foot" of $N(\tau)$ $k$ times around a curve $K''$ on
$\Sigma^*$  which is parallel to $K$ on $\Sigma$.
 
\QED

\vskip20pt

\hskip30pt\BoxedEPSF{Fig.6.1.eps scaled810}
\bigskip
\centerline{Fig. $6.1$}
\bigskip

\tag\Sameforlinks
\state Remark. \rm Notice that  Lemma \Propositionstab \ remains correct if we replace the knot
$K$  by a $q$-component link $L$ and \lq\lq $1$-isotopic" by \lq\lq$q$-isotopic". This is
because one can do the same operations as explained in the last proof, with one  stabilization 
required for each component of $L$ .
\bigskip

\tag\peelrem
\state Remark. \rm   Let $K \subset S^3$ be a knot carried by some generalized
trellis  $T$, and let $\Sigma$ be the associated  trellis Heegaard surface. Then $T$  
is a spine of the handlebody  $H_1 = N(T)$.   After  drilling out a properly embedded 
arc $\tau \subset H_1$ and isotopying the boundary of the new handlebody slightly 
off $K$  as in the proof of Lemma 6.1(a), we can enlarge the tunnel $N(\tau)$  and 
thus \lq\lq peel off" $N(K)$  from $H_1$ to get a handlebody $H'_1$.   The boundary
$\Sigma' = \partial H'_1$ is isotopic in $S^3 - K$ to the Heegaard surface $\Sigma^*$ from 
the proof of Lemma 6.1.  Compared to $H_1$ the new handlebody  $H'_1$  contains 
an extra handle,  namely the neighborhood of the  peeled off knot  $K$. The core $K$  
of this extra  handle is connected to $T$ by a small arc $\sigma$ which we call the 
{\it stem} of  the knot, see Fig. $6.2$ below.

\vskip20pt

\hskip30pt\BoxedEPSF{Fig.6.2.eps scaled 800} 
\bigskip
\centerline{Fig. 6.2}
\bigskip

For the next proof we need to introduce a new operation on the trellis $T$,
called a {\it top} (or {\it bottom})  {\it horizontal edge slide}.  It consists of
taking the top (or the bottom) vertex $w$ of a vertical edge $e$ which is the
outer-most  (say left-most) vertex on some horizontal line of $T$, and sliding $w$
along that horizontal line to the other end. The edge $e$ is isotoped into a position 
behind the original trellis $T$, and its top (bottom)  endpoint is now the right-most 
endpoint of the  new horizontal line.  A picture is  given in Fig. 6.3.

 \vskip20pt

\hskip40pt\BoxedEPSF{Fig.6.3.eps scaled770}
\bigskip
\centerline{Fig. 6.3}
\bigskip

Notice that whenever $K$ is carried by $T$,  then such an edge horizontal slide
will induce an isotopy of $K$, as we keep $K$ on $\partial N(T)$ throughout the 
edge slide. With respect to the new {\it twisted trellis}, obtained after the horizontal 
edge slide, the long arc of $K$ on the back of the old trellis (at the hight level of the
horizontal line of $T$ along which the horizontal edge slide was performed) has now 
become a short right-most horizontal arc on the front  of the twisted trellis, while the
former left-most short horizontal arc on the front has now become a long horizontal 
arc on the back.

\tag\sweepout
\state Lemma. Let $K$ be a knot carried by a generalized trellis $T$, and let $v$ 
be the left-most vertex of the top horizontal line of $T$.  Then there exists a finite 
sequence of horizontal edge slides on $T$ such that the resulting twisted trellis
$T'$ has the following property:  \hskip 50pt 
\indent If one starts to slide a point $x \in K$ close to $v$ along the subarc of $K$ which 
runs parallel to the whole length  of the top  horizontal line of $T'$ and then once around 
$K$, then, for every horizontal  layer of $T'$, $x$ will completely traverse the top and the 
bottom horizontal line, before it ever crosses over more than one vertical edge (called a ``special
edge" from that horizontal layer.

\give Proof. The point  $x$  starts moving along  $K$  by traveling first along all
of the top horizontal line of $T$, and then down, winding around the right-most
vertical edge $e$ of the first horizontal layer. According to whether the 
twist coefficient of $e$ is odd or even, the point $x$ has to continue by sliding towards 
the left or towards the right.  Correspondingly we apply a top horizontal edge slide 
to all vertical edges of the second layer which have their top endpoints to the  left 
(or to the right) of the bottom vertex of the special edge $e$, and similarly a bottom horizontal 
edge slide to all edges of the first layer which have their bottom endpoints to 
the left (or to the right) side of the bottom vertex of $e$.

As a consequence $x$ ends up on the subarc of $K$ which completely traverses
the second horizontal line of $K$, and we have to consider the possibility that the
endpoint of this horizontal line is the bottom endpoint of a vertical edge in the
first layer.  In this case $x$ will move again up into this first layer, until it
eventually reaches a vertex on the second horizontal line which contains the top
endpoint of a vertical edge $e'$ from the second horizontal line.  Then  $x$ slides 
down on $K$ into the second layer, winding around the special edge $e'$, and we have 
to repeat  the procedure just explained, with $e'$ replacing $e$.  This is repeated finitely 
many times until we have swept out all horizontal lines of $K$.

\QED

\tag\Propositionsequi
\state Proposition.  For any knot  $K \subset S^3$, carried by a generalized trellis
$T$, the  associated trellis Heegaard surface $\Sigma$ is $1$-isotopic in
$\Sigma({1\over k}) = K({1 + k \, a(K) \over k})$,  for any $k \in \integers$, to a multiple 
stabilization of the canonical top or bottom Heegaard surfaces $\Sigma_{top}$ or $\Sigma_{bot}$
associated to $K$.
\bigskip

\give Proof.   We first change the trellis  $T$  (and the knot $K$ accordingly)
by doing horizontal edge slides so that it satisfies the conclusion of Lemma \sweepout\/. 
Let $\Sigma$ be the trellis Heegaard surface for the pair  $(S^3,K)$ given by the 
resulting twisted trellis, still called $T$, and let $K'$ be the core of the surgery filling 
of $\Sigma({1\over k})$. We stabilize $\Sigma$ in $\Sigma({1\over k})$ to get 
the vertical Heegaard surface $\Sigma^*$ of $\Sigma({1\over k}) - \inter N(K') =
S^3 - \inter N(K)$ as in Lemma \Propositionstab \ (a).  By Lemma \Propositionstab \ (b) 
this is isotopic in $S^3 - N(K)$ to the Heegaard surface $\Sigma'$ defined by $H'_1$ in 
Remark \peelrem\/. Let $\sigma$  be the stem as defined there.  We will prove the
proposition by describing a sequence of slides of the edges of $T$ and of $\sigma$.
We always think of $\Sigma'$ as of the boundary of a small regular neighborhood of the 
handlebody spine which is isotoped along throughout the sequence of slides.

We first introduce a slide of the stem $\sigma$ in $S^3$ which keeps  the knot endpoint 
of  $\sigma$ on $K$ and the trellis endpoint on $T$.  This can be done in  such a way  that the
stem is always a short straight arc, for example by keeping it throughout the slide perpendicular
to the edge of the trellis along which the trellis endpoint of $\sigma$ is moving.  In particular 
this shows that we can freely choose the starting position of  $\sigma$. We choose as
starting vertex for the trellis endpoint of $\sigma$ the top left corner vertex $v$ of $T$, and for
the knot endpoint the point $x$ given by Lemma \sweepout. 

The {\it stem slide} is now defined by sliding $\sigma$  in the described fashion so that its knot
endpoint moves around
$K$ exactly once. Note that, by the time it comes back to  $v$,  every edge  $e$ of 
$T$  has been  traversed precisely twice by the trellis endpoint of $\sigma$.  

Next we introduce, for every edge  $e$ of $T$, a {\it second comming slide} as follows:  Immediately 
after traversing $e$ for the second time, i.e., with the stem positioned at the "second" endpoint
$x$ of $e$, we interrupt the above stem slide. We isotope the edge $e$ of the trellis along the
knot, by sliding its endpoint $x$ first over the  stem $\sigma$ and then back along $K$,  so that 
$e$ is now  replaced by a new stem which is attached to the other endpoint of the former edge
$e$. As this is done after the second  and last time the trellis endpoint visits the edge $e$, we can
complete  the above stem slide of $\sigma$ once around $K$, although the edge $e$  is now missing in 
$T$.    

We define the {\it second coming procedure} as follows: Perform the stem slide, but as the stem
$\sigma$ slides around $K$ do the second comming slide to every edge $e$ of $T$. This
creates lots of new stems and eliminates eventually all edges of the trellis.
We now investigate more precisely the effect of this second coming procedure on the 
edges of the trellis:

Before doing this procedure, vertical edges of the trellis had either $2, 3$ or $4$  adjacent
horizontal edges, depending on their location in the trellis.  Consider a vertical edge $e$ which
had $4$ adjacent horizontal edges. Notice that it follows from the horizontal edge slides performed
at the beginning of the proof in accordance with Lemma \sweepout\/ that in the stem slide the trellis 
endpoint of $\sigma$ traverses each of the $4$ horizontal edges at least once before it crosses over the 
vertical edge $e$ for the first time . Hence the first passage of the trellis endpoint through $e$ will
produce precisely one stem at one of the endpoints of $e$, and none at the other.  The second
passage through $e$ will produce a stem at each of the two endpoints of $e$.  Thus for every
vertical edge $e$ with $4$ adjacent horizontal edges the second coming procedure gives precisely a 
single stem at one of the endpoints of $e$ and a double stem at the other.  Note that this double stem 
has none of its endpoints on the string of the knot $K$ which runs horizontally on the back of
$\Sigma$.  Instead, it connects the two strings of $K$ which wind around the vertical edge
$e$.

By similar considerations the same conclusion holds for vertical edges with $3$ adjacent horizontal
edges, if the analogous assumption is satisfied.  This includes the horizontal edge in the first layer with
$v$ as top vertex (even if it has only $2$ adjacent horizontal edges), as can be seen directly from the 
fact that the original stem, in final position, will be placed with trellis endpoint at $v$.

Observe now that, as a consequence of the horizontal edge slides performed on $T$ and $K$ 
at the beginning of this proof, the only vertical edges in $T$ with only $2$ adjacent horizontal 
edges are possibly the special edges from Lemma \sweepout\/ or the edge with endpoint $v$.  
Hence, by the time the trellis endpoint of $\sigma$ has  returned to the starting vertex $v$, there 
will  be a double stem at the top or the bottom  of the corresponding twist box for all vertical 
edges except for the one special  edge in every layer.

The Heegaard surface $\Sigma''$ isotopic to $\Sigma'$ which results from the second coming
procedure is hence obtained from $K$ by introducing the tunnel system given by all the  double
stems (the single ones can be deleted without changing the isotopy class of $\Sigma''$).
Thus it follows from Lemma \horizontalarcs\/  that $\Sigma''$ is isotopic to a multiple
stabilization of either, the top or the bottom  vertical Heegaard surface
$\Sigma_{top}$ or $\Sigma_{bot}$.

\QED

We can now apply the above proposition to knots with flypes and obtain:
\bigskip

\tag\Stabflype
\state Theorem. Let $T$ be a generalized trellis with an interior pair of edges, and
let $T(n)$ be the trellis obtained from an $n$-flype at this edge pair. Let $K \subset S^3$ be a 
knot carried by $T$.  Then for all $n \in \integers$ the trellis Heegaard surface 
$\Sigma(n) \subset K({1 + k \, a(K) \over k})$ 
of genus $g(T(n)) = g(T) + 2n$,  given by the trellis $T(n)$,  is $1$-isotopic
to a multiple stabilization of the canonical top or bottom Heegaard surfaces $\Sigma_{top}$ or
$\Sigma_{bot}$ of $K({1 + k \, a(K) \over k})$, defined by the trellis projection of 
$K$ before the flypes.
 
\give Proof.  This is an immediate consequence of the last Proposition 6.4. and of Proposition
\Flypestab .

\QED

\tag\Ericthesedgwick
\state Remark.  \rm An alternative proof of the last theorem can be given by combining
a result of Sedgwick [Se] with Proposition \Propositionsequi .  Notice that Sedgwick's proof
applies to a more general situation than the one given by trellisses, since it is a local proof. 
Consequently, it is not possible to deduce the statement of Proposition \Propositionsequi\/ by his
methods, as that statement is of global nature.
\smallskip
\bigskip
\noindent{\bf Proof of Theorem 0.3.}  The theorem follows directly from 
Proposition \Propositionsequi \/ and Theorem \Stabflype.

\QED

\noindent{\section {Stabilizing canonical vertical Heegaard splittings}
\bigskip
\bigskip

In this section we investigate the question of how many stabilizations are necessary so that the 
canonical top and bottom  Heegaard splittings of a knot  $K \subset S^3$, given as  a 2n-plat, 
become isotopic. It was proved by Hagiwara (see [Ha]) that $n-1$ stabilizations always suffice. 
We give a new proof of this result and show that in many cases one can do with considerably 
fewer stabilizations. The following notion has been introduced, with minor technical variation, 
in [LM].

\bigskip

\tag\widthdef
\state Definition. \rm (a)  A  $2n$-braid $b$ will  be said to have {\it total width} $r \in  \naturals$ 
if in its standard projection $\pi (b) \subset P$ (obtained from  $b$ by  replacing every crossing  
by a node), every monotonically descending  path in $\pi (b)$ connecting the $i$-th  strand at the top 
to the $k$-th strand on the bottom satisfies $i - k \leq 2r -1$ , and  $r $ is the smallest
such number.\hfil \hskip 400pt
\noindent(b) A $2n$-plat projection of a knot $K$ will be said to have {\it total width} 
$r \in  \naturals$ if the underlying $2n$-braid has total width $r$ .
\bigskip

Clearly for any $2n$-plat one always has $0\leq r \leq n-1$.
If $r = 0$ then the braid in question defines the $n$-component
unlink.  If $m$ is the number of horizontal layers of the plat, then one has $r \leq (m + 1) / 2$.

   We prove:

\tag\hagprop
\state Proposition. For every knot or link $K \subset S^3$ of width $r$ the
two canonical Heegaard splittings of $S^3 - \inter N(K)$ are $r$-isotopic.

\give Proof. For each $i$ with $r +1 \leq i \leq n-1$ we consider the tunnel
system consisting of the tunnels
$\tau_1 , ..., \tau_{i} , \rho_{i-r} , \rho_{i-r+1}  , ...,\rho_{n-1}$.
We claim that this system is isotopic to the system
$\tau_1 , ..., \tau_{i-1} , \rho_{i-r} ,\rho_{i-r+1}  , ...,\rho_{n-1} , \eta$,
where $\eta$ is a trivial tunnel. Assuming this claim it follows from  the
symmetry between the top and from bottom tunnels that the system
\centerline{$\tau_1 , ..., \tau_{i-1} , \rho_{i-r} , \rho_{i-r+1},...,\rho_{n-1} ,\eta$}
 is isotopic to the system

\centerline{$\tau_1 , ..., \tau_{i-1} , \rho_{i-r-1} , \rho_{i-r},...,\rho_{n-1} $.}

\noindent Thus we conclude inductively that 

\centerline{$\tau_1 , ..., \tau_{n-1} , \rho_{n-r-1} , \rho_{n-r},...,\rho_{n-1} $}

\noindent is isotopic to

\centerline{$\tau_1 , ..., \tau_{r} , \rho_{1} , \rho_{2} ,...,\rho_{n-1} $.}

\bigskip

\noindent  It follows from Lemma 5.6 that these systems are just the canonical systems with
$ r$ trivial tunnels added.

It remains to prove the above claim.  The assumption on the total width of the $2n$-plat $K$
implies that there is no monotonocally decending path connecting the  $2i$-th strand on 
the top to the $(2(i-r) - 2)$-strand on the bottom of the plat.  In other words, the left most
monotonically descending path  $\gamma$ in $\pi (K)$  which starts at the top of the 
$2i$-th  string must end at the bottom of some $k$-th string with $k \geq 2(i - r)-1$.  We
consider the handlebody  $W = N(K \cup \tau_1  \cup ... \cup  \tau_{i-1}  \cup
\rho_{i-r}  \cup\rho_{i-r+1}   \cup  ... \cup \rho_{n-1})$ and,  for all  $j = i-r, ..., n-1$,
we introduce cocore disks $D_j$ for the tunnels $\rho_j$ (see Fig. 5.1).

Consider now an equatorial $2$-sphere $S$ which intersects  $K$ just below
the top bridges and cuts off a $3$-ball $B$ with  $n$ unknotted arcs 
$t_1, \dots , t_n$ (the top bridges), as indicated in  Fig. 7.1.

\vskip15pt

\hskip50pt\BoxedEPSF{Fig.7.1.eps scaled700}\bigskip
\medskip
\centerline{Fig. 7.1}
\bigskip

Let $\beta$ be an arc in  $S$
which is isotopic relative boundary to the top tunnel $\tau_i$ . Consider an isotopy  $\cal I$ of  $S$
determined by moving the sphere monotonically down, so that at each level we have a horizontal
$2$-sphere, to a level just above the the bottom bridges. The isotopy  $\cal I$ moves the
intersection points  $\cup _{j = 1}^n \{t_j  \cap S \}$ in such a way on  $S$ that it
braids the arcs $t_j$ according to the strands of the given $2n$-plat $K$.  Let $\beta'$ be
the image of $\beta$ after the isotopy  $\cal I$.

We can assume without loss of generality that each crossing of $\pi (K)$  lies
on a distinct height level, called a \lq\lq critical" level. The left-most descending 
path $\gamma\,$ in $\pi (b)$, defined above, determines at each horizontal
level a split of $S$ into a \lq\lq left" and \lq\lq right" half. (To be precise, $S - \{\infty\}$
is split along $\pi ^{-1}(\gamma)$, where $\pi : \reals^3 \to P$  is the orthogonal projection.)
As we move $S$ by the isotopy  $\cal I$ through a critical level we see iteratively that  $\cal I$  can be
chosen so that the arc $\beta$ is alway contained in an $\epsilon$-neighborhood of the
right half of  $S$  determined by $\gamma\,$, were $\epsilon$ is smaller than the distance between any
two strands of the plat.  In particular, when  $S$  has reached the bottom level, then
the obtained arc $\beta'$  is positioned entirely to the right of the $(2(i-r) - 2)$-th strand. 
Thus we can isotope $\beta'$ on  $\partial W$ to become a small trivial arc $\eta$ by 
sliding it across  the cocore disks $D_{i-r} , ..., D_{n-1}$ of the tunnels 
$\rho_{i-r},  \dots, \rho_{n-1}$. This proves the claim and finishes the proof of the proposition.

\QED

\section References.
\bigskip
\bigskip

[BO]  M. Boileau, J.P. Otal; {\it Scindements de Heegaard et groupe des homeotopies

\hskip 24pt des petites varietes de Seifert},  Invent. Math. 106,  (1991) 85 - 107

[BZ1] M. Boileau, H. Zieschang; {\it Heegaard genus of closed orientable
Seifert

\hskip 28pt 3 - manifolds},  Invent. Math. 76,  (1984) 455 - 468

[BZ2] G. Burde, H. Zieschang; {\it Knots} , De Gruyter studies in Mathematics 5,

\hskip 28pt (1985)

[CG] A. Casson, C. McA. Gordon; {\it  Reducing Heegaard splitting}, Top. and
App. 27,

\hskip 26pt (1987)  275 - 273

[Ha] \hskip 3pt Y. Hagiwara; {\it  Reidemeister - Singer distance for unknotting
tunnels of a knot},

\hskip 26pt  Kobe J. of Math. 11,  (1994)  89 - 100

[He]  \hskip 3pt J. Hemple; {\it 3 - manifolds}, Ann. Math. Studies 86, (1976)
Princeton University

\hskip 26pt Press

[Ka]  \hskip 3pt A. Kawauchi; {\it Classification of pretzel knots},   Kobe J. of
Math.  2, (1985) 1 - 22

[Ko]  \hskip 3pt T. Kobayashi; {\it  A construction of  3-manifolds whose
homeomorphism classes of

\hskip 26pt   Heegaard splittings have polynomial growth}, Osaka J. Math.  29, (1992)  \hfil\break 
\indent \hskip 26pt   653 - 674

[LM] M. Lustig, Y. Moriah; {\it Closed incompressible surfaces in complements
of wide

\hskip 26pt knots and links}, to appear in Topology and its Applications.

[LM1] M. Lustig, Y. Moriah; {\it Nielsen equivalence for Fuchsian groups and and Seifert 

\hskip 28pt fibered spaces}, Topology 30, (1991)  191 - 204

[LM2] M. Lustig, Y. Moriah;{\it Generalized Montesinos Knots, Tunnels and $\cal N$-Torsion}, 

\hskip 30pt   Math. Ann. 295,   (1993)  167 - 189  

[Ly] \hskip 4pt  H. Lyon; {\it Knots without unknotted incompressible spanning surfaces}, Proc. Amer.

\hskip 24pt  Math. Soc. 35, (1972)  617 -  620

 [MS] Y. Moriah, J. Schultens;  {\it Irreducible Heegaard splittings of Seifert
fibered spaces

\hskip 24pt are either  vertical or horizontal}, to appear in Topology.

 [MR] Y. Moriah, H. Rubinstein; {\it Heegaard structure of negatively curved
3-manifolds},

\hskip 26pt   to appear in Comm. in Ann. and Geom. .

[Ro] \hskip 3pt D. Rolfsen; {\it Knots and links}, Math. Lect. Series 7, (1976)
Publish or Perish

[RS] \hskip 2pt  H. Rubinstein, M. Scharlemann; {\it Comparing  Heegaard
splittings of non-Haken

\hskip 24pt   3-manifolds}, preprint.
\eject
 [Se]  \hskip 3pt E. Sedgwick; {\it An infinite collection of Heegaard spittings that are equivalent

\hskip 24pt after one stabiliztion} Math. Ann.  308, (1997)  65 - 72

 [Sc]  \hskip 3ptJ. Schultens;{\it The classification of Heegaard splittings for
(closed  orientable

\hskip 22pt surface) $\times S^1$}, Proc. London Math. Soc. 67, (1993)  401 - 487

[ST] M. Scharlemann, A. Thompson, {\it Heegaard splittings of (surfaces)$\times
I$  are standard,}

\hskip 21pt Math. Ann. 295,  (1993)  549 - 564

[St]  \hskip 3pt M. Stocking, {\it Almost normal surfaces in 3 - manifolds},
preprint.

\bigskip

\bigskip
\bigskip

\obeyspaces Martin Lustig                                            Yoav
Moriah

\obeyspaces Institut fur Mathematik                               Department of
Mathematics

Rhur-Universitat Bochum                             Technion, Haifa  32000,

Postfach 10 21 48                                        Israel

D-4630 Bochum 1                                       ymoriah@techunix.technion.ac.il

Germany

\end


\bigskip
\tag\trellisthm
\state Theorem 0.2. Let $T$ be a generalized trellis of genus $g$ and let $K = K(A) \subset S^3$ be a 
knot or link carried  by $T$ with twist matrix $A$. Assume that all coefficients $a_{i,j}$ of $A$  which
are different from $\infty$ satisfy  $|a_{i,j}| \geq 3$ and that there is an interior pair of edges $(e_{i,j},
e_{i,j +1})$ of  $T$ with twist coefficients $|a_{i,j }|, |a_{i,j  + 1}| \geq 4$.   Then for all $n \in
\integers$ the trellis $T(n)$ obtained from an $n$-flype at this edge pair defines  a trellis Heegaard 
splitting for the pair $(S^3, L)$ which is horizontal and of genus $g + 2n$. In particular, if $K$ is a 
knot, then for all the manifolds $\Sigma({1\over k})$ obtained from $1\over k$-Dehn filling at $K$, with  
$|k| \geq 6$, this induces a strongly irreducible Heegaard splitting of genus $g+2n$.